\documentclass[12pt, twoside]{amsart}

\swapnumbers\newtheorem{theorem}{Theorem}[section]

\newtheorem{cory}[theorem]{Corollary}
\newtheorem{prop}[theorem]{Proposition}
\newtheorem{lemma}[theorem]{Lemma}

\newtheorem{remarks}[theorem]{Remarks}

\numberwithin{equation}{section}

\usepackage{amssymb,latexsym, amsmath, amscd}

\def \sirc{{\raise0.2ex \hbox{$\scriptstyle \circ$}}}

\def \wt{\widetilde}
\def\pa{\partial}

\def\po{\parindent 0pt}
\def\p{{\po \it \s Proof. }}

\def\s{\smallskip}
\def\m{\medskip}

\def\IM{\operatorname {Im}}

\def\wh{\operatorname{Wh}}
\def\wp{\operatorname{Wh}(\pi)}
\def\grad{\operatorname{grad}}

\def\cat{\operatorname{cat}}

\def\RR{\mathbb {R}}

\def\ZZ{\mathbb{Z}}
\def\zp{\ZZ[\pi]}

\def\ga{\alpha}

\def\eps{\varepsilon}

\def\gf{\varphi}

\def\gl{\lambda}

\def\gt{\tau}

\begin{document}

\title[Functions on $h$-cobordisms]{On the minimal number of critical points
of functions on $h$-cobordisms}

%\date{\today}

\author{P. E. Pushkar}
\address{P. E. Pushkar, Math. Department University of Toronto, } 
\email{ppushkar@math.toronto.edu}
\author{Yu. B. Rudyak}
\address{Yu.\ B.\ Rudyak, Math. Inst. Univ. Heidelberg, Im Neuenheimer Feld
288, 69120 Heidelberg, Germany} \email{rudyak@mathi.uni-heidelberg.de}

\begin{abstract} Let $(W,M_0,M_1)$ be a non-trivial $h$-cobordism
(i.e., the Whitehead torsion of $(W,M_0)$ is non-zero) with $\dim W >5$. We
prove that every smooth function $f: W  \to [0,1]$, $f(M_0)=0, f(M_1)=1$ has at
least 2 critical points. This estimate is sharp: $W$ possesses a function as
above with precisely two critical points.
\end{abstract}

\maketitle

\section*{Introduction}

Let $(W,M_0,M_1)$ be an $h$-cobordism, \cite{Mi}. Here $W$ is always
assumed to be smooth, connected and compact and $M_i, i=0,1$ is always
assumed to be closed. Recall that an $h$-cobordism $(W,M_0,M_1)$ is called
trivial if there is a diffeomorphism $(W,M_0,M_1)\cong (M\times [0,1], M_0,
M_0)$. We say that a function(not necessarily Morse) $f: W \to [0,1]$ is {\it
regular} if  $f^{-1}(M_0)=0$, $f^{-1}(M_1)=1$ and both values 0 and 1 are
regular values of $f$. It is well known that an $h$-cobordism $(W,M_0,M_1)$ is
trivial if and only if $W$ possesses a regular function without critical
points. In this note we prove the following theorem.

\m
{\parindent 0pt \bf Theorem}
{\it Let $(W,M_0,M_1)$ be a non-trivial $h$-cobordism with $\dim W \ge 6$.
Then every regular function on $W$ has at least two critical points.
Moreover, this estimate is sharp: $W$ possesses a regular function with
precisely two critical points.} 

\m
We denote by $I$ the closed interval $[0,1]$.

\m
{\parindent 0pt \bf Acknowledgments.} This work was partially supported by the
Fields Institute for Research in Mathematical Science, Toronto, the
Simplectic Topology, Geometry and Gauge Theory Program. The first author was
also supported by RFBR grant 99-01-01109 NWO grant 047-008-005, NSERC PD
Fellowship.

\section{Preliminaries}\label{prel}

 Let $f: W \to I$ be a regular Morse function on an $h$-cobordism $(W,
M_0,M_1)$. Choose a Riemannian metric on $W$ and consider integral 
trajectories for the vector field $-\grad f$, the so-called anti-gradient
trajectories. We say that an anti-gradient trajectory $a=a(t)$ is {\it a
special trajectory from $p$ to $q$} if $\lim_{t \to -\infty} a(t)=p$ and
$\lim_{t \to +\infty} a(t)=q$ where $p$ and $q$ are critical points of $f$
such that the index of $p$ is one more than the index of $q$. We can and shall
assume that the number of special trajectories is finite (this is true for
generic function and metric).

For every critical point of $f$ we fix orientations of unstable disks
(left-hand disks in terminolgy of \cite{Mi}). Then every unstable sphere (in a
certain level) gets an orientation. Moreover, every stable sphere gets a
coorientation, i.e., an orientation of its normal bundle in the corresponding
level set. Now, for every special trajectory $a$ from $p$ to $q$ we define the
number $\eps(a)=\pm 1$ as follows. Take $c \in ]f(q), f(p[$. Then our
trajectory $a$ meets the level $f^{-1}(c)$ in a certain point $x$, which is a
point of transversal intersection of the corresponding stable and unstable
spheres. We define $\eps (a)$ to be the intersection index at $x$.

\section{Whitehead torsion} 

Given a ring $R$, we define a based $R$-module to be a free finite generated
left $R$-module $M$ with a fixed $R$-free basis.

\m Recall the definition of the Whitehead torsion of an $h$-cobordism
$(W,M_0,M_1)$. Given a group $\pi$, let $A=A(\pi)$ denote the set of long
exact sequences  
$$ 
\CD \cdots @>>> C_n @>\pa_n>> C_{n-1}
@>>> \cdots @>>> 0 
\endCD 
$$
such that each $C_i$ is a based $\zp$-module and all but finite number of
modules $C_i$ are zero modules. Furthermore, each $\pa_i$ is a $\zp$-module
homomorphism. Let us call the exact sequence of based $\zp$-modules trivial if
it has only two non-zero terms and the corresponding isomorphism is given by
the identity matrix. The term-wise direct sum operation converts $A$ into an
abelian semigroup. Let $R$ be the equivalence relation on $A$ generated by the
following operations:

\begin{itemize}
\item interchanging of the elements;
\item replacement of a basis element by the sum of this element with the
multiple of another basis element; 
\item addition of the
trivial exact sequence;  
\item multiplication of any basis element by the
element $\pm g, g \in \pi$. 
\end{itemize}

\m The above mentioned operation in $A$ induces a group structure in
$A/R$. This groups is called the{\it Whitehead group} of $\pi$ and is denoted
by $\wp$, \cite{Mil}. It turns out to be that $\wp$ is a functor of $\pi$. In
particular, every homomorphism $\gf: \pi \to G$ induces a homomorphism
$\wh(\gf): \wp \to \wh(G)$. Namely, the homomorphism $\gf$ yields the
homomomorphism $\ZZ[\gf]: \zp \to \ZZ[G]$ of group rings, and for every based
$\zp$-module $C$ we can form the based $\ZZ[G]$-module
$C\otimes_{\ZZ[\gf]}\ZZ[G]$. The sequence $\{C_n
\otimes_{\ZZ[\gf]}\ZZ[G]\}$ turns out to be exact because all the $C_n$'s are
free, etc.

\m For every $h$-cobordism $(W,M_0,M_1)$ with $\pi_1(W)=\pi$ the Whitehead
torsion $\gt(W,M_0,M_1)\in \wp$ is defined as follows. Consider a regular
Morse function $f: W \to I$, Riemannian metric, etc. as in \S \ref{prel}  . Fix
a point $x_0\in W$ and, for every critical point $p$ of $f$, choose a path
$u(p)$ from $x_0$ to $p$. Given a special trajectory $a$ from $p$ to $q$, we
define a path $v=v_a: I \to W$ as follows. Let $\gl(t): \RR \to ]0,1[$ be a
function such that 
$$ 
\lim \gl(t)_{t \to 0}=-\infty, \quad \lim \gl(t)_{t \to
1}=+\infty. 
$$ 
We set $v(0)=p$, $v(1)=q$, $v(t)=a(\gl (t))$. Now, consider the
loop $u(p) \sirc v \sirc (u(q))^{-1}$ (where $\sirc$ denotes the composition
of paths) and define $g(a)\in \pi =\pi_1(W)$ as the based homotopy class of
the loop.   

\m Let $p_1, \ldots, p_k$ be all the critical points of the index $n$. Define
$C_n$ to be the free $\zp$-module generated by symbols $[p_1], \ldots, [p_k]$.
In other words, $C_n$ consists of formal linear combinations  
$$
\sum_{i=1}^k\ga_i[p_i], \quad \ga_i\in \zp. 
$$
We define the differential $\pa_n: C_n \to C_{n-1}$ to be a $\zp$-module homomorphism such that 
$$
\pa_n[p]=\sum_q\sum_{a\in T(p,q)}\eps(a)g(a)[q]
$$
where $q$ runs over all critical points of the index $n-1$ and $T(p,q)$ is
the set os special trajectories from $p$ to $q$.

\m It follows from the Morse theory that $H_*(\{C_n, \pa_n\})= H_*(\wt W, \wt
M_0)$ where $((\wt W, \wt M_0)$ is the universal covering of the pair
$(W,M_0)$. Since $M_0$ is a deformation retract of $W$, we conclude that $\wt
M_0$ is a deformation retract of $\wt W$, and therefore the complex $\{C_n,
\pa_n\}$ is acyclic, i.e. the sequence 
$$ 
\CD \cdots @>>> C_n @>\pa_n>> C_{n-1}
@>>> \cdots @>>> 0 
\endCD 
$$
is exact. Thus, the above sequence determines a certain element $\gt = \gt
(W,M_0)\in \wp$, the so-called {\it Whitehead torsion} of the $h$-cobordism
$(W,M_0,M_1)$. 

\m According to well-known Barden--Mazur--Stallings Theorem, \cite{K,M,S}, an
$h$-cobordsim $(W,M_0,M_1)$ with $\dim W \ge 6$ is trivial if and only if
$\tau(W,M_0)=0$.

\begin{lemma}\label{triv}
 Suppose that an $h$-cobordism $(W,M_0,M_1)$ possesses a regular Morse
function $f$ such that all the critical points and special trajectories of $f$
are contained in a simply connected domain $U$ of $W$. Then $\tau(W,M_0)=0$
\end{lemma}

\p Since $\tau(W,M_0)$ does not depend on the choice of the based point $x_0$
and the paths $u(p)$, we can assume that $x_0\in U$ and every path $u(p)$
belongs to $U$. Then, for every special trajectory $a$, $ga$ is the neutral
element of $\pi=\pi_1(W)$. Thus, $$ \tau(W,M_0)\in \IM\{\wh(j): \wh\{e\} \to
\wp\} $$ where $j:\{e\} \to \pi$ is the inclusion of the trivial subgroup. But
it follows from the elementary linear algebra that $\wh\{e\}=0$, see e.g.
\cite{Mil}. Thus, $\tau(W,M_0)=0$.

\section{Proof of the theorem}

 Let $f:M \to \RR$ be a smooth function (not necessarily Morse) on a
Riemannian manifold $M$. Let $U$ be an open ball in $M$ and suppose that $U$
contains precisely one critical point $o$. 

\begin{lemma}\label{close}
There exists a regular function $g$ which is $C^{\infty}$-closed to $f$ in the
Whitney topology and such that every special $g$-trajectory is contained in $U$
whenever its ends are contained in $U$.   
\end{lemma}

\p Let $D(r)=\{m\in M\bigm| d(m,o) < r)$ where $d$ is the distance function
on $M$. We can and shall assume that the injectivity radius at $p$ is at least
one and that $U=D(1)$. Then there constants positive constants $C$ and $E$
such that, for every function $g$ which is $C^{\infty}$-closed to $f$, the
following estimates holds in $D(1)\setminus D(1/2)$: 
$$ 
|\grad g|\ge E, \quad
|L_{\grad f}\,d(m,o)|\le C. 
$$  
Choose a function $g$ closed to $f$ let $p$ and
$q$ be two critical points of $g$ which belongs to $U$. Suppose that there is
a special trajectory $a(t)$ from $p$ to $q$ which meets the boundary of
$D(3/4)$. We claim that in this case 
$$ 
g(p)-g(q)\ge \frac{E^2} {4C}. 
$$
Indeed, since $L_{\grad f}\,d(m,o)|\le C$, we conclude that  
$$
a\left[t-\frac 1 {4C}, t+\frac 1 {4C}\right]
$$
does not meet $D(1/2)$ whenever $a(t)\notin D(3/4)$ . So, if $a(t_0)\notin
D(3/4)$ then %
\[
g(p)-g(q)\ge \int_{t_0-\frac 1 {4C}}^{t_0+\frac 1 {4C}}dg(a(t)) 
=\int_{t_0-\frac 1 {4C}}^{t_0+\frac 1 {4C}}|\grad g|^2dt \ge \frac{E^2} {4C}.
\]

Now we can finish the proof as follows. Since $f$ has only one critical
point, there exists $g$ closed to $f$ and such that $g(p)-g(q)$ is small enough
for every critical point $p$ and $q$ of $g$. This is a contradiction. 
\qed 

\begin{cory}\label{one}
If an $h$-cobordism $(W,M_0,M_1)$ possesses a regular function $f$ with one
critical point $p$, then $\tau(W,M_0)=0$. In particular, if $\dim W\ge 6$ then
the $h$-cobordism is trivial. 
\end{cory}

\p Because of Lemma \ref{close}, we can perturb the function $f$ in a small
neighborhood of the critical point and get a function $f_1$ such that all its
critical points and special trajectories belong to a disk neighborhood of $p$.
Now the result follows from Lemma \ref{triv}. 
\qed

\begin{prop}\label{takens}
Every $h$-cobordism $(W,M_0,M_1)$ , $\dim W\ge 6$ possesses a regular
function with at most $2$ critical points. 
\end{prop}

\p Consider a regular Morse function $f: W \to I$. Asserting as in
\cite[Lemme 1]{K} and \cite[\S 4]{Mi}, we can modify $f$ and to get a regular
Morse function which has at most two critical levels $a,b$, $a<b$ and index of
each of critical points is equal to 2 or 3. Because of this, every critical
level is path connected. Now, following \cite[Th. 2.7 and Prop.2.9]{T}, we can
contract the critical points in each of levels and get a regular function with
at most 2 critical points. 
\qed

\m
Clearly, Corollary \ref{one} and Proposition \ref{takens} together imply the Theorem.

\begin{remarks}\rm 1. Asserting as in \ref{one}, one can show that, for every
regular function $f$ on a non-trivial $h$-cobordism, the number of critical
levels of $f$ is at least 2 provided that all the critical points of $f$ are
isolated.

\m 2. Every $h$-cobordism $(W,M_0,M_1)$ possesses a regular function with 1
critical level. Namely, choose collars of the boundary components and define 
$f$ to be constant on complements of collars and depending on ``vertical''
coordinate only for collars. In greater detail, consider a smooth function %
\begin{equation*} \gf: I \to I, \quad \gf(t)=\left\{ \begin{array}{lcl}
&t  &\text{ if }\ 0\le t \le\eps/4  \text{ or }\ 1-\eps/4 \le t \le 1, \\
&1/2  &\text { if }\ \eps/2 \le t \le 1-\eps/2
\end{array}
\right.
\end{equation*}
for $\eps>0$ small enough. Choose collars $M_0 \times [0,\eps]$ and
$M_1\times [1-\eps,1]$ and define $f: W \to I$ by setting 
\begin{equation*}
f(x)=\left\{
\begin{array}{lcl}
&\gf(t) &\text{ if}\ x=(m,t)\in M_0\times [0,\eps],\\
&\gf(t) &\text{ if}\ x=(m,t)\in M_1\times [1,1-\eps],\\
&1/2  &\text { else. } 
\end{array}
\right.
\end{equation*}

\m 3. Every trivial $h$-cobordism ($M\times I,M,M)$ possesses a regular
function with 1 critical point. Indeed, consider a function $\gf: M \to I$
such that $\gf^{-1}(1)$ is a point $m_0$ (and therefore $m_0$ is a critical
point of $\gf$) and define $$ f: M \times I \to I, \quad
f(m,t)=(t-1/2)(1-\gf(m))+\gf(m)(t-1/2)^3. $$  It is easy to see that $f$ has
just one critical point $(m_0,1/2)$.

\m 4. Notice that, for every $h$-cobordism $(W,M_0,M_1)$, the relative
Lusternik--Schnirelmann category $\cat(W,M_0)=0$, while every regular
function on any non-trivial $h$-cobordism $(W,M_0,M_1)$ has at least two
critical points.   

\m 5. It is easy to see that, because of the collar theorem, the regularity 
condition for $f$ in the Theorem can be weaken as follows: $f(M_0)=0$ and
$f(M_1)=1$. 
 
\end{remarks}

\end{document}